\documentclass[12pt]{amsart}

\usepackage{graphicx}
\usepackage{amsfonts}
\usepackage{hyperref}
\usepackage{textcomp}
\usepackage{amssymb}

\makeatletter
\def\url@leostyle{%
  \@ifundefined{selectfont}{\def\UrlFont{\sf}}{%
    \def\UrlFont{\footnotesize\ttfamily}}\Url@do
}
\makeatother

\usepackage{graphicx}

\newcommand{\union}{\cup}

\def\N{\mathcal{N}}
\def\P{\mathcal{P}}

\def\D{\mathcal{D}}

\renewcommand{\v}{{\mathbf v}}

\usepackage{amsmath}
\newtheorem{thm}{Theorem}

\newtheorem{lem}[thm]{Lemma}

\newtheorem{prop}[thm]{Proposition}

\newtheorem{obs}[thm]{Observation}

\newcommand{\Wnor}{{W_{\mathsf{nor}}}}

\newcommand{\Wmis}{{W_{\mathsf{mis}}}}

\date{}

\title{Monotonic Sequence Games}
\author[Albert et al.]{M.~H.~Albert}
\address{Department of Computer Science \\%
University of Otago} 
\email{malbert@cs.otago.ac.nz}
\author[]{R.~E.~L.~Aldred}
\address{Department of Mathematics and Statistics \\%
University of Otago} 
\author[]{M.~D.~Atkinson}
\address{Department of Computer Science \\%
University of Otago} 
\author[]{C.~C.~Handley}
\address{Department of Computer Science \\%
University of Otago} 
\author[]{D.~A.~Holton}
\address{Department of Mathematics and Statistics \\%
University of Otago} 
\author[]{D.~J.~McCaughan}
\address{Department of Mathematics and Statistics \\%
University of Otago} 
\author[]{B.~E.~Sagan}
\address{Department of Mathematics\\%
Michigan State University}

\begin{document}
\begin{abstract}
In a monotonic sequence game, two players alternately choose elements of a sequence from some fixed ordered set. The game ends when the resulting sequence contains either an ascending subsequence of length $a$ or a descending one of length $d$. We investigate the behaviour of this game when played on finite linear orders or $\mathbb Q$ and provide some general observations for play on arbitrary ordered sets.
\end{abstract}

\maketitle


\section{Introduction}

Monotonic sequence games were introduced by Harary, Sagan and West in \cite{HSW:83}. We paraphrase the description of the rules as follows:
\begin{quote}
From a deck of cards labelled with the integers from $1$ through $n$, two players take turns choosing a card and adding it to the right hand end of a row of cards. The game ends when there is a subsequence of $a$ cards in the row whose values form an ascending sequence, or of $d$ cards whose values form a descending sequence.
\end{quote}
The parameters $a$, $d$, and $n$ are set before the game begins. There are two possible methods for determining the winner of the game. In the {\em normal}\/ form of the game, the winner is the player who places the last card (which forms an ascending or descending sequence of the required length). In the {\em mis{\`e}re}\/ form of the game, that player is the loser. In \cite{HSW:83} these are called the {\em achievement}\/ and {\em avoidance}\/ forms of the game respectively.

As a consequence of the Erd{\H{o}}s-Szekeres theorem \cite{ES:35}, the game cannot end in a draw if $n > (a-1)(d-1)$. It is therefore natural to attempt to classify the parameters $(a,d,n)$ according to whether the first player can force a win, the second player can force a win, or either player can ensure at least a draw. Some results towards such a classification were presented in \cite{HSW:83} and the problem of extending and generalising these results was posed there and by Sagan in \cite{Sagan:01}.

In this paper we will report on some progress on this and related problems. As regards the original game, we have been able to extend the computer-assisted analysis to decide many instances which were left open in \cite{HSW:83}. We also provide some general results concerning the long run behaviour of these games (that is, for fixed $a$ and $d$ but large $n$). However, most of the work reported here deals with variations of the original game. In particular, we consider the case where the deck of cards is $\mathbb Q$ rather than a finite linear order. Finally, we examine some other variations of the game obtained either by relaxing the rules, or by playing with a deck of cards that is partially ordered. We list some open problems in the final section of the paper. 

We adopt, and in some cases adapt, the notation and terminology of {\em Winning Ways}\/ \cite{BCG1:04, BCG2:04, BCG3:04,BCG4:04} in discussing our results.   This differs somewhat from that used in \cite{HSW:83} so, where necessary, we will also provide translations of the results from that paper.  

\section{The general framework}

Any version of the monotonic sequence game specifies at the outside, a {\em deck}\/ $D$ which is simply some partially ordered set, and two positive integer parameters $a$ and $d$ which we call the {\em critical lengths}\/ of ascending and descending sequences respectively. There are two players, $A$ and $B$ (for convenience in assigning pronouns, $A$ is assumed to be male and $B$ female), who alternately choose an element which has not previously
 been chosen from the deck and add it to a sequence whose elements consist of the cards chosen up to this point. This sequence will be called the {\em board}. Conventionally, $A$ plays first while $B$ plays second. In the {\em basic}\/ form of the game
the board is constructed from left to right. That is, if the current board is $b  c \cdots v$ and the next player chooses a value $w \in D$ then the new board is $b c \cdots v w$. An ascending subsequence of length $a$ or a descending subsequence of length $d$ of the board is called a {\em critical sequence}. As soon as the board contains a critical sequence the game ends.
In {\em normal}\/ play, the winner is the player whose move terminated the game. In {\em mis{\`e}re}\/ play that player is the loser. We henceforth assume that $a,d \geq 2$ since the cases $a = 1$ or $d = 1$ are completely trivial.
If the deck is exhausted without creating a critical sequence, then the game is considered drawn. If the deck is infinite then the game is also considered drawn if play proceeds without termination. By default we assume that normal play is being considered unless otherwise noted. 

\begin{prop}
If $D$ is finite and contains a chain of length greater than $(a-1)(d-1)$, or $D$ is infinite and contains no infinite antichain then no draws are possible in either normal or mis{\`e}re play.
\end{prop}

\begin{proof}
In the first case any supposedly drawn board would contain all the elements of the specified chain. However, by the Erd{\H{o}}s-Szekeres theorem any such sequence contains a critical sequence. In the latter case a similar result follows from the well-known observation that, as a consequence of Ramsey's theorem, any infinite sequence of elements from a partially ordered set contains an infinite subsequence which is either ascending, descending, or an antichain. Since the last possibility is ruled out by hypothesis, one of the former two must apply, and the play producing that sequence could not have been drawn.
\end{proof}

\begin{obs}
If $D$ has a fixed-point-free order-preserving involution then the second player can force at least a draw.
\end{obs}

$B$'s strategy is to play the image of $A$'s move under the involution, {\em unless}\/ she has an immediate win available. Since no chain can involve both a point and its image she thereby never plays a suicidal move, that is one which makes it possible for Alexander to win the game on his next turn, and hence she cannot lose.

\begin{obs}
If $a = d$ and $D$ has a fixed point free order reversing involution $i$ with the property that whenever $x$ and $x^i$ are comparable, one is minimal and the other maximal, then the second player can force at least a draw.
\end{obs}

Again the strategy for $B$ is to play a winning move if one exists, and otherwise the image of $A$'s previous move. The minimality/maximality criterion guarantees that in the resulting sequence of plays no chain can arise using both $x$ and $x^i$ unless $a = d = 2$ which is trivially a second player win.

This observation applies to play on the cube $2^n$ or equivalently on  the lattice of subsets of a set. In particular it is easy to check that for $a = d= 3$ play on $2^3$ is a second player win though cooperatively the two players can play to a draw. 

Since $D$, $a$ and $d$ are fixed parameters of any particular game, all the relevant information about a position is contained in its board. A board which could arise in play may not have a proper prefix containing either an ascending sequence of length $a$ or a descending sequence of length $d$. Subject to this condition we may define the {\em type}\/ of a board to be one of $\N$, $\P$ or $\D$. We say that the type is $\N$ (next) if the player whose turn it is to move (that is, the next player) has a winning strategy. The type is $\P$ (previous) if the previous player (that is, the player who is not next) has a winning strategy. Finally, the type is $\D$ (drawn) if each player has a strategy that guarantees her or him at least a draw.

A board which contains the entire deck or which contains a critical sequence is called a {\em terminal board}. A terminal board containing a critical sequence is of type $\P$ in normal play and $\N$ in mis{\`e}re play, while a terminal board that does not contain a critical sequence is of type $\D$. Otherwise, the type of a board, $X$, is determined by the set of types of the boards that can be obtained in one further move. We call these boards the {\em children}\/ of $X$. If this set contains any board of type $\P$ then the type of $X$ is $\N$. If {\em all}\/ the boards in this set are of type $\N$ then the type of $X$ is $\P$. Otherwise, the type of $X$ is $\D$.

These rules may not be immediately sufficient for determining the type of an arbitrary board when arbitrarily long plays or even draws with infinite play are possible. However, even in this case the boards are partitioned into the three types above. The algorithm for performing the partitioning is to begin by labelling all the terminal boards according to the winning conditions. Then inductively any currently unlabelled boards which either have a child of type $\P$, or all of whose children have type $\N$, are labelled appropriately. After completing this induction, any boards remaining unlabelled are of type $\D$.

Our principal goal will be to determine the type of the empty board -- that is, to determine whether the first player has a winning strategy, or failing that, whether he can force a draw. We denote this type by $\Wnor(a,d,D)$ for normal play, or $\Wmis(a,d,D)$ for mis{\`e}re play.

\section{Double bumping}

Given a sequence of distinct elements $\v = v_1 v_2 \cdots$ from a linearly ordered set $C$ there is a well known algorithm due to Schensted \cite{Schensted:CJM61} which determines (explicitly) the length of the longest increasing subsequence of any prefix $v_1 v_2 \cdots v_k$ and (implicitly) the elements of such a sequence. This is sometimes called the ``bumping'' algorithm. An increasing  sequence ${\mathbf w} = w_1 w_2 \cdots w_m$ is maintained as the
elements of $\v$ are processed in order. When $v_i$ is processed, $\mathbf w$ is modified as follows: if $w_m < v_i$ then $v_i$ is appended to $\mathbf w$; otherwise $v_i$ bumps (that is, replaces) the smallest element of $w$ that is larger than $v_i$.

It is easy to check that, after processing $v_1 v_2 \cdots v_k$ the element $w_j$ of $\mathbf w$ is the least maximum element of an ascending subsequence of $v_1 v_2 \cdots v_k$ of length $j$.
In particular, the length of $\mathbf w$ is equal to the length of the longest ascending subsequence obtained to that time.

Of course there is a dual algorithm that allows one to keep track of the length of the longest descending subsequence. In this version an element is either prepended to the sequence being maintained (if smaller than all the elements of the sequence), or it bumps the immediately smaller element. 

For the purposes of analysing some forms of the monotonic sequence game it will be useful to be able to combine these two algorithms into a single one. However, in doing so, we need to keep track of whether the elements in the single ordered sequence which we are maintaining represent elements of the ascending or descending type -- that is, whether an element takes part in the sequence $\mathbf{w}$ of the original algorithm, the corresponding sequence $\mathbf{m}$ in the dual algorithm, or both.

Initially we will do this by marking the elements with overlines (if they belong to $\mathbf{w}$), underlines (if they belong to $\mathbf{m}$) or both (if both). Thus we maintain a single marked sequence which we shall call the {\em recording sequence}. The double bumping form of the combined algorithm can then be described as follows.
\begin{itemize}
\item
Initially set the recording sequence to be empty, and process the elements of the permutation in order from left to right.
\item
Repeatedly, until the permutation is exhausted:
\begin{itemize}
\item
insert the next element of the permutation into the recording sequence with both an underline and an overline (maintaining the increasing order of the recording sequence); 
\item
delete the first overline if any to its right and the first underline if any to its left; 
\item
remove any naked elements (ones which no longer have an underline or an overline). 
\end{itemize}
\end{itemize}

For example, when we process the permutation $514263$ in this way we obtain:
\newcommand{\ol}[1]{{\overline{#1}}}
\newcommand{\ul}[1]{{\underline{#1}}}
\newcommand{\oul}[1]{{\ol{\ul{#1}}}}

\[
\oul{5} \: \rightarrow \: 
\oul{1} \ul{5} \: \rightarrow \: 
\ol{1} \oul{4} \ul{5} \: \rightarrow \: 
\ol{1} \oul{2} \ul{4} \ul{5} \: \rightarrow \: 
\ol{1} \oul{2} \ul{4} \oul{6} \: \rightarrow \: 
\ol{1} \ol{2} \oul{3} \ul{4} \ul{6}
\]

Frequently the precise identity of the elements of the recording sequences will not be important, but only their type (that is, what decoration they have). This remark will be exactly true when we deal with monotonic sequence games on $\mathbb Q$, and is still of some relevance in the case of monotonic sequence games on finite chains. For typographical purposes it is easier to record type sequences as colours rather than bars, and so we will also call them {\em colour sequences}. Specifically we associate the colour Blue with an underline, and Red with an overline. Elements having both underlines and overlines will be called Purple. An element is reddish if it is Red or Purple, and bluish if it is Blue or Purple. The process of the double bumping algorithm on the permutation above, purely in terms of colours is:
\[
P \: \rightarrow \: 
PB \: \rightarrow \: 
RPB \: \rightarrow \: 
RPBB \: \rightarrow \: 
RPBP \: \rightarrow \: 
RRPBB.
\]
Of course the length of the colour sequence corresponding to a permutation is not more than the length of the permutation itself. Different permutations can easily have the same colour sequence (e.g. $231$ and $213$ both have colour sequence $PP$) and indeed permutations of different sizes can have the same colour sequence (e.g. $312$ and $2143$ both have colour sequence $RPB$).

It is clear that only some sequences of colours can occur as a result of applying the double bumping algorithm. We call such colour sequences {\em admissible}. In terms of colour, when we add a new element, we insert a Purple somewhere in the sequence and remove the red tinge from the first reddish element to the right (deleting it entirely if it were Red) and the blue tinge from the first bluish element to the left. In particular, a colour sequence can never begin with Blue nor end with Red. In fact we can completely characterise the admissible colour sequences. Recall that a {\em factor}\/ of a sequence is a subword consisting of a block of consecutive elements from the sequence.

\begin{prop}
The language of  admissible colour sequences consists precisely of the empty sequence, together with those sequences which contain at least one $P$, do not begin with $B$ nor end with $R$, and do not contain $RB$ as a factor.
\end{prop}

\begin{proof}
Necessity is relatively straightforward. Each insertion leaves a $P$ so a non-empty admissible sequence must contain a $P$. Of the remaining conditions, the first two conditions are obviously preserved by any legitimate insertion. To see that the final condition is preserved as well consider an insertion which supposedly creates an $RB$ factor. It could not create both the $R$ and the $B$ since only an insertion between those two elements could do that. Suppose, without loss of generality, that the newly created element was the $B$. Then previously that element was represented by a $P$. But in order to eliminate its reddish tinge, the insertion would have had to be after any preceding $R$, so we could not get the $RB$ factor as claimed.

The proof of sufficiency is by induction. We show that if $w$ is a non-empty sequence of the form described, then there is some parent word $v$ also of the form described such that $w$ can be obtained from $v$ by the bumping algorithm. That this suffices is based on the observation that for {\em any}\/ starting word $u$ (admissible or not), after $ad+1$ bumps the resulting word must contain at least $d+1$ bluish or $a+1$ reddish (red or purple elements). Thus the backwards chain of parents from $w$ is bounded in length by the product of the number of bluish elements and the number of reddish elements in $w$, and can only terminate in the empty sequence which is admissible.

If $w = P$ the result is clear, so we may assume that the length of $w$ is at least two. Suppose first that  $w = P u$. If $u = B u'$ let $v = P u'$ (which still has the form required) and note that $v$ produces $w$ by an insertion on the left hand side. If $u$ begins with a $P$ or an $R$ let 
$v = R u$ which is admissible and produces $w$ by an insertion on the left hand side.

Now suppose that $w = R^i P u$ with $i > 0$ and let $w' = Pu$. Then $w'$ is admissible, and by the case just proven we can find $v'$ which produces $w'$ by an insertion into the first position. Let $v = R^i v'$. Then $v$ produces $w$ by insertion after the first block of $R$'s.
\end{proof}

The number of non-empty admissible words is enumerated by the sequence of alternate Fibonacci numbers:
\[
1, \, 3, \, 8, \, 21, \, 55, \, 144, \ldots
\]
This is easily established by standard transfer matrix approaches or by the observation that the association:
\[
R \rightarrow 01 \quad B \rightarrow 10 \quad P \rightarrow 00
\]
almost provides a bijection between admissible colour sequences and binary sequences of even length which contain no consecutive $1$'s.

\section{Finite chains}

In this section we assume throughout that the deck is a finite chain which, for convenience, we take to be:
\[
[ n ] = \{ 1, 2, \ldots, n \}
\]
with the usual ordering. This was the basic situation investigated by Harary, Sagan and West in \cite{HSW:83}. On the theoretical front we have relatively little to add to their results in this area, however, we have extended their computational results considerably.

\begin{prop}
For fixed $a$ and $d$ both sequences $\Wnor(a,d,[n])$ and $\Wmis(a,d,[n])$ for $n = 1, 2, 3, \ldots$ are eventually constant.
\end{prop}

\begin{proof}
Since we know that a play of the game with parameters $a$ and $d$ cannot last more than $(a-1)(d-1) + 1$ moves, the existence of a winning strategy for either player, in either termination condition, can be expressed as a first order sentence in the language of linear orders. Consider, for example, the case of a first player win in normal play. In this case this sentence begins with an existential quantifier, followed by a long alternation of quantifiers representing the moves which might be chosen by the two players. These quantifiers are followed by a quantifier free formula expressing the condition ``the first ascending sequence of length $a$ or descending sequence of length $d$ arising in this play occurred after a move made by the first player''. The other cases are all similar.

However, it is well known that the theory of finite linear orders admits quantifier elimination (see 
\cite{Hodges:93}, specifically sections 2.7 and A.6 and their exercises).  In particular, any sentence in this language is either true in $([n], <)$ for all sufficiently large $n$ or false in $([n], <)$ for all sufficiently large $n$. Since one of the statements ``the game is of type $\N$'' and ``the game is of type $\P$'' must be true for every $n > (a-1)(d-1)$, it must be the case that the same one is true for all sufficiently large $n$.
\end{proof}

The proof above is a little unsatisfying from the standpoint of attempting to understand the structure of the monotonic sequence game played with a finite deck. By essentially recreating the quantifier elimination for the theory of finite linear orders but tailoring it to the situation at hand we can make it somewhat more concrete. As a side-effect we obtain improved bounds for the onset of the ``long term behaviour'' of such games.

Specifically, consider boards that arise in the play of the monotonic sequence game. Suppose that the colour sequence at this point $c_1 c_2 \cdots c_k$. There is an associated sequence of gap lengths $g_0, g_1, \cdots g_k$ where $g_i$ is the number of cards remaining in the deck between the elements representing $c_{i-1}$ and $c_i$. Note that this is not necessarily the same as the difference between these elements minus one, as some of the
 intervening elements may have been played earlier but no longer form part of the colour sequence. 

The basic idea of the argument is to divide gaps into two categories {\em large}\/ and {\em small}. All gaps whose length is larger than a certain number (which may depend on the colour sequence and the position of the gap relative to that sequence) will be considered large. We aim to show that if two boards have the same colour sequences and corresponding gaps are either both large or both small and of equal length then we can emulate the following play in one game within the other game and vice versa. This {\em Tweedledum-Tweedledee argument}\/ then establishes that the two games have the same outcome type (and in fact the same nim-value or Grundy number). 
The first part of the argument must establish just what the bounds are for large gaps.

Imagine for the moment that the next play of the game will be a card from the deck that lies in some particular gap. Among the values in the board below this card there will be some maximal increasing sequence whose length, $r$, is the number of reddish elements lying below the gap. Likewise there is some maximal decreasing sequence on the board whose length, $b$ is the number of bluish elements lying above the gap. Within this particular gap, the game will certainly end if we create an increasing sequence of length $a-r$ or a decreasing one of length $d-b$. That is, within the gap we are essentially playing a game with parameters $a-r$ and $d-b$ (the play within this gap may influence plays in other gaps, but only by reducing their associated parameters). Suppose that we temporarily let $B(x,y)$ denote some value which is ``big enough'' to define a large gap for parameters $x$ and $y$. A play into such a gap leaves two gaps, a lower one with parameters $x$ and $y-1$ and an upper one with parameters $x-1$ and $y$. Since we must ensure that we can match small gaps exactly and create corresponding large gaps it will be sufficient to have:
\[
B(x,y) \geq B(x-1, y) + B(x, y-1) + 1.
\]
If we choose equality and note that we may take $B(x,1) = B(1,y) = 1$ then simple algebraic manipulation shows that we may choose:
\[
B(x,y) = 2  {{x+y-2} \choose {x-1}} -1.
\]
Henceforth we take this as the definition of $B(x,y)$ and hence of what constitutes a large gap.

\begin{prop}
\label{prop:equivalence}
For fixed $a$ and $d$, any two boards having the same colour sequence with the property that corresponding gaps are either both large, or otherwise equal have the same outcome type.
\end{prop}

\begin{proof}
As promised, the proof is what is known as a ``Tweedledum-Tweedledee'' argument (in the language of \cite{BCG1:04}) or a ``back and forth'' argument (in the language of \cite{Hodges:93}). The idea is that any move made in either position has one or more matching moves on the other position which preserve the equality of colour sequences and corresponding gaps. Specifically, a move in a small gap is mirrored by the obvious corresponding move of the other position. A move in a large gap leaves either large gaps on either side or one small gap and one large gap. In either case there is a corresponding move in the other position leaving two large gaps, or one small gap (of the same size) and a large gap.

Suppose, for the sake of argument, that the first position has a second player winning strategy. We devise a second player winning strategy in the second position as follows. Given a move in the second game to which we must reply, we consider a matching move in the first game. Our strategy there will dictate a certain response to this move. We make the matching response in the second game. Proceeding in this way, we cannot fail to win in the second game (in fact we will win in precisely the same number of moves as we win the matched sequence of plays in the first game). All the other cases are very similar.
\end{proof}

In particular, any two games beginning with an empty board and having decks of size $2 {{a+d-2} \choose {a-2}}-1$ or larger must have the same outcome type. As indicated by the computations below, this bound appears to be somewhat extravagant, though not as much so as the na\"{i}ve bound arising from a direct translation of the quantifier elimination for the theory of finite linear orders which would be $2^{(a-1)(d-1)}-1$. 

\subsection{Computational results: normal play}

We will assume throughout that $a \geq d$ because the outcome type for parameters $(a,d)$ is the same as that for parameters $(d,a)$. We begin by recapitulating results from \cite{HSW:83} recast into our notation.

If $d = 2$ then any move other than the smallest remaining element at that time gives your opponent a ``win in one''. So the outcome type is determined by the parity of $a$ and we have:
\[
\Wnor(a,2,n) = 
\left\{
\begin{array}{ll}
\D & \mbox{if $n < a$} \\
\N & \mbox{if $a \leq n$ is odd} \\
\P & \mbox{if $a \leq n$ is even.}
\end{array}
\right.
\]

For $d = 3$ then, depending on parity, the first player can choose to play either the largest or second largest element of the deck as his first move. This more or less reduces the game to the $d = 2$ case, and  provided that $n > a$ and $a$ is even, or $n > a+1$ and $a$ is odd $\Wnor(a,3,n) = \N$, with the remaining cases being drawn. 

Finally, \cite{HSW:83} showed that $\Wnor(4,4,n) = \N$ for $n \geq 9$. A winning strategy is to play near the middle, and to ensure after your second move that all remaining moves must be the smallest or largest remaining element.

We implemented a straightforward game tree traversal algorithm to determine the outcome type of the empty board for various combinations of the parameters $(a,d,n)$. Although the observations made in the proof of Proposition \ref{prop:equivalence} could improve the efficiency of this algorithm (by storage and re-use of previously computed outcomes for equal or equivalent colour and gap sequences) such time improvement would come at significant cost in space, and complexity of the underlying code. Since we could extend the results of \cite{HSW:83} considerably using just the raw improvement in computing power between 1983 and now, we did not choose to pursue these improvements. Our program permitted computations with deck sizes up to 20 in a few minutes on a standard desktop machine. Note that whenever a type $\P$ position is found, two other positions are immediately known to be of type $\N$, namely:
\begin{eqnarray*}
\lefteqn{\Wnor(a,d,n) = \P \: \Rightarrow} \\
& & \mbox{$\Wnor(a+1,d,n+1) = \N$ and $\Wnor(a,d+1,n+1) = \N$}
\end{eqnarray*}
since the first player can reduce the game to the preceding case by playing the smallest (respectively largest) element as his first move.

We give our new computational results in the following form: first we specify the smallest non-drawn game of that type and its winner; then a sequence of values until we (appear) to reach an eventually constant block. Thus, the first line below means that for $n \leq 10$ $\Wnor(5,4,n) = \D$, and for $11 \leq n \leq 20$ $\Wnor(5,4,n) = \N$.
\[
\begin{array}{c}
(5,4,11) \in \N \\
(6,4,14), (6,4,15) \in \P, (6,4,16) \in \N \\
(5,5,15) \in \N \\
(7,4,15), (7,4,16) \in \N, (7,4,17) \in \P, (7,4,18) \in \N.
\end{array}
\]

\subsection{Computational results: mis{\`e}re play}

We also computed results for mis{\`e}re play. In this case it appears to be true that the game is drawn much less frequently, and so the data include some more interesting observations. In this case, the table below lists the sequence of outcome types for the various combinations of parameters $a$ and $d$ with $n$ ranging from $1$ through $20$.
\[
\begin{array}{cc|llll}
a & d & \multicolumn{4}{c}{\mbox{Mis{\`e}re winner}}  \\ \hline
3 & 3 & \D\D\D\N\N & \N\N\N\N\N & \N\N\N\N\N & \N\N\N\N\N \\
4 & 3 & \D\D\D\D\P & \N\N\N\N\N & \N\N\N\N\N & \N\N\N\N\N \\
5 & 3 & \D\D\D\D\D & \N\P\P\P\N & \N\N\N\N\N & \N\N\N\N\N \\
6 & 3 & \D\D\D\D\D & \D\D\N\N\N & \N\N\N\N\N & \N\N\N\N\N \\
7 & 3 & \D\D\D\D\D & \D\D\D\P\N & \N\N\N\N\N & \N\N\N\N\N \\
8 & 3 & \D\D\D\D\D & \D\D\D\D\N & \P\P\N\N\N & \N\N\N\N\N \\
9 & 3 & \D\D\D\D\D & \D\D\D\D\D & \D\N\N\N\N & \N\N\N\N\N \\
4 & 4 & \D\D\D\D\D & \N\P\N\N\N & \N\N\N\N\N & \N\N\N\N\N \\
5 & 4 & \D\D\D\D\D & \D\D\N\N\N & \P\N\N\N\N & \N\N\N\N\N \\
6 & 4 & \D\D\D\D\D & \D\D\D\D\N & \P\N\P\P\N & \N\N\N\N\N \\
7 & 4 & \D\D\D\D\D & \D\D\D\D\N & \P\N\N\N\N & \N\N\N\N\N \\
8 & 4 & \D\D\D\D\D & \D\D\D\D\D & \D\N\N\N\D & \N\N\N\N\P \\
9 & 4 & \D\D\D\D\D & \D\D\D\D\D & \D\D\D\N\P & \N\P\N\N\D \\
5 & 5 & \D\D\D\D\D & \D\D\D\D\D & \D\N\P\N\P & \N\N\N\N\N \\
6 & 5 & \D\D\D\D\D & \D\D\D\D\D & \D\N\P\N\N & \N\P\N\N\N \\
7 & 5 & \D\D\D\D\D & \D\D\D\D\D & \D\N\D\N\P & \N\P\N\N\P
\end{array}
\]

Most of the blocks of trailing $\N$'s do seem to represent long run behaviour. The evidence supporting this is that the smallest winning first move is also constant across these blocks.

The $a = 8,9$, $d = 4$ cases seem particularly interesting. First of all, with $a = 8$ there is the interposed $\D$ at $n = 15$. Thus, with a 14 or 16 card deck the first player can force the second player to make an
 ascending sequence of size 8 or a descending one of size 4 but with a 15 card deck he cannot! A further oddity of this sequence concerns
 the fact that for $a = 9$, $d = 4$, the {\em second}\/ player wins $n = 15$. This means that the second player can force the first to create an ascending sequence of length 9 or a descending one
 of length 4 in a 15 card deck, but can't force an ascending sequence of length 8 or a descending sequence of length 4 in the same deck. Why can't the first player simply follow an ``at least draw'' strategy from the latter case to get the same result in the former case? Because there is a hidden assumption in this strategy -- that the second player will never create an ascending sequence of length 8 or a descending sequence of length 4 either.

\subsection{Computation: further remarks}

As noted above the program used to obtain these results was exceedingly straightforward. Essentially, every response to every move was examined from lowest to highest. Only when a response of type $\P$ was found (permitting the current board to be labelled as $\N$) was any pruning done. Likewise, no heuristic choices of responses were considered. This alone would probably improve the efficiency of the program considerably since it was observed that in many cases if $y$ was a good countermove to first move $x$ (and was quite different from $x$) then it was also a good countermove to $x+1$. Secondly, storage and re-use of previously computed results, or  some form of ``orderly'' generation based on Proposition \ref{prop:equivalence} would permit even more pruning. For example, the first three moves $10, 5, 20$ and $10, 20, 5$ result in identical colour sequences and gaps, so have the same outcome type.

However, beyond some obvious observations and conjectures which we propose in the final section, our opinion is that the data (particularly for the mis{\`e}re version) suggest rather ``noisy'' behaviour for small values of $n$. So, the benefits of pursuing these optimisations seems rather limited.

\section{Dense linear order}

We now consider playing the monotonic sequence game with $\mathbb Q$ (or any other dense linear order without endpoints) as the deck. 

\begin{prop}
For any $a, d \geq 1$, $\Wnor(a,d,{\mathbb Q}) = \Wmis(a-1, d-1, {\mathbb Q})$.
\end{prop}

\begin{proof}
In order to win the normal game, you cannot ever create an ascending chain of length $a-1$ or a descending chain of length $d-1$ since your opponent would then have the opportunity to win immediately. Conversely, if your opponent creates such a sequence on the board then you can win immediately. So the outcome of the mis{\`e}re $(a-1, d-1, {\mathbb Q})$ game is the same as that of the normal $(a,d,{\mathbb Q})$.
\end{proof}

We note that the proposition above requires only that the deck not have a maximal or minimal element. Owing to this proposition we restrict our attention to the normal form of the game.

The outcome type of a particular board depends only on the relative ordering among the elements currently on the board. This is clear, since with two boards having the same relative ordering among their elements, there is an order preserving bijection from $\mathbb Q$ to itself which maps one board to the other. Any strategy which applies to the first board, then also applies to the second by taking its image under this bijection. However, in fact all that we need to know in order to determine the outcome of a game is the colour sequence of the board. As noted previously, different boards and even boards of different sizes can have the same colour sequence.

\begin{prop}
For any $a, d \geq 1$, and playing with $\mathbb Q$ as a deck, the outcome type of a particular board is determined by the colour sequence of that board.
\end{prop}

\begin{proof}
First we observe that the colour sequence of the board is sufficient to determine whether or not the game has ended since the length of the longest ascending (descending) sequence on the board is equal to the number of reddish (bluish) elements of its colour sequence. 

Next we note that given the colour sequence of a board, the possible colour sequences which can be obtained by making a single move are determined. Any move involves the insertion of a $P$ somewhere in the existing colour sequence, and then ``first higher red reduction'' and ``first lower blue reduction''. Moreover, because the deck is dense, any such insertion {\em can}\/ be made. 

So, in terms of determining the outcome, we need only know the colour sequence of the current board, exactly as claimed.
\end{proof}

In considering the basic form of the monotonic sequence game with parameters $(a,d,{\mathbb Q})$ we will work almost exclusively with the colour sequences. We define the {\em children}\/ of a colour sequence to be all those sequences that can be obtained from it in a single move. A colour sequence is {\em terminal}\/ if it contains $a$ reddish, or $d$ bluish elements.

As before, we will assume that $a \geq d$ and for a few values of $d$ we are able to determine the type of the general game with parameters $(a,d,{\mathbb Q})$.

\begin{thm}
For $d \leq 5$ the types of the monotonic  sequence games with parameters $(a,d,{\mathbb Q})$ are as follows:
\begin{enumerate}
\item
For $a \geq 2$, $\Wnor(a,2,{\mathbb Q}) = \P$.
\item
For $a \geq 3$,  $\Wnor(a,3,{\mathbb Q}) = \N$ precisely when $a$ is odd.
\item
For $a \geq 4$, $\Wnor(a,4,{\mathbb Q}) = \N$.
\item
For $a \geq 5$, $\Wnor(a,5,{\mathbb Q}) = \N$.
\end{enumerate}
\end{thm}

\begin{proof}
Throughout the argument we consider an equivalent version of the monotonic sequence game with parameters $(a,d,{\mathbb Q})$. In this version, a suicidal move i.e. one which creates an ascending sequence of length $a-1$ or a descending sequence of length $d-1$ on the board is forbidden, unless forced. Since the player with a winning strategy in the original game will never make a suicidal move, and the other player may choose not to so until forced, the outcome type of the modified game is the same as that of the original.

For parameters  $(a,2,{\mathbb Q})$ the game is truly trivial, since the very first move is suicidal.

For the parameters $(a,3,{\mathbb Q})$, any move below an element already played is suicidal. So, in the modified form, the two players alternately add to an increasing sequence, and clearly the first player wins only if $a$ is odd.

Now consider parameter sequences of the form $(a,4,{\mathbb Q})$. We will show that the set of colour sequences representing non-terminal $\P$-positions in this game is: 
\[
P_4=  \{P, R^{a-3} P B\} \cup 
\{ R^i P^2 \, | \,  0 \leq i \leq a-5 \} .
\]
To establish this result we must show that for any position which can arise in the play of $(a,4,{\mathbb Q})$, if it is not in $P_4$ then it has a child which is in $P_4$ or a terminal position, and if it is in $P_4$ then there is no such child. The second part is easily checked.

Suppose that we have a colour sequence $w$ which is not terminal and not in $P_4$. If it has three or more bluish elements, then it has a terminal child. Suppose that $w$ has exactly one bluish element. Then it is of the form $R^i P$ for some $0 < i \leq a-2$. If $i = a-2$ it has a terminal child. If $i < a-2$ then an insertion just before the last $R$ yields $R^{i-1} P^2$ which is in $P_4$ unless $i = a-3$. In that case, inserting before the $P$ yields $R^{a-3} P B$. 

Next consider the case where $w$ has two bluish elements, both purple. Ignoring positions with terminal children, it must be of the form $R^i P R^j P$ where either $j > 0$ or $i = a-4$. If $j > 0$ inserting before the last $R$ yields $R^{i+j} P^2$ while inserting before the last $P$ yields $R^{i+j+1} P B$ and one of these two is in $P_4$. If $j = 0$ and $i = a-4$ then inserting between the two $P$'s yields $R^{a-3} P B$. 

Finally consider the case of one purple and one blue element. Then $w$ is  $R^i P B$ for some $i$. If $i \leq a-5$ then moving at the right hand end produces $R^i P^2$, while if $i = a-4$, moving just after the $P$ produces $R^{a-3} P B$.

Thus for the parameters $(a,4,{\mathbb Q})$ we have established that $P \in \P$ and hence the initial position is in $\N$.

We give a similar argument for the parameter sequences of the form $(a,5,{\mathbb Q})$. In this case though we do not provide an exhaustive listing of the type $\P$ non-terminal colour sequences, but only a sufficient set of these. By this we mean that we provide a set $P_5$ of colour sequences, and an argument that the following conditions hold:
\begin{itemize}
\item
$P \in P_5$;
\item
if $w \in P_5$ and $v$ is a child of $w$, then $v$ has a child which is either terminal or  in $P_5$;
\item
no $w \in P_5$ has a terminal child.
\end{itemize}
This establishes that $P \in \P$, since from any position not in $P_5$ the player whose turn it is to move can simply take either the immediate win, or the move guaranteed by the second of the conditions above.
We take 
\begin{eqnarray*}
P_5 &=& \left\{P, \, RPB \right\} \union \\
{} && \left\{R^i PRPB \, : \, 0 \leq i \leq a-6 \right\} \union \\
{} && \left\{R^i RPBP \, : \, 0 \leq i \leq a-6 \right\} \union \\
{} && \left\{R^{a-5} P^3, \, R^{a-3} P B^2\right\}.
\end{eqnarray*}
The first part of the verification is routine. From the initial position $P$ the second player can ensure that after her play, the resulting code will be $RPB$ by always replying ``in the second position''. Likewise from $RPB$ she can always guarantee that her opponent's next move will be from one of $PRPB$ or $RPBP$.

Now suppose that $0 \leq i \leq a-7$ and that a single move has been made from $R^i  PRPB$ or $R^i  RPBP$. If this move occurs below the first $P$ it creates a descent of length $4$ and can be countered by an immediate win (i.e. it is suicidal). In all of the remaining cases there is a counter move to one of $R^{i+1} PRPB$ or $R^{i+1} RPBP$.

If a single move has been made from $R^{a-6} PRPB$ or $R^{a-6} RPBP$ which is non-suicidal, then again there are only a few positions near the end of the colour sequence that need to be examined, and each of these allows a response to either $R^{a-5} P^3$ or $R^{a-3} PB^2$.

The final cases to consider are moves from $R^{a-5} P^3$. There are only two non-suicidal moves and they both permit replies to $R^{a-3} P B^2$.

\end{proof}

We have strong experimental evidence that the monotonic sequence game with parameters $(a,d, {\mathbb Q})$ and $a,d \geq 4$ {\em always}\/ has type $\N$. Computation has established this result for $4 \leq d \leq 8$ and any $a$ with $d \leq a \leq 16$. We can establish this result for the symmetrical form of the game:

\begin{thm}
Let $a \geq 4$. The monotonic sequence game with parameters $(a,a,{\mathbb Q})$ has type $\N$.
\end{thm}

\begin{proof}
The argument we provide uses a form of {\em strategy stealing}\/ together with symmetry. That is, we show that if the second player had a winning strategy then the first player could appropriate it for his own use. This contradiction implies that it must be the first player who has a winning strategy.

If the result were false then the type of the colour sequence $P$ would have to be $\N$. As the moves from $P$ to $RP$ and $PB$ are symmetrical (under order reversal) both these positions would have to be of type $\P$. 

In particular the two children $PP$ and $RPB$ of $RP$ would both lie in $\N$. The children of $PP$ are $PBP$, $RPB$ and $PRP$. By assumption, $RPB \in \N$. By symmetry $PBP$ and $PRP$ have the same type, so these two positions would have type $\P$. The children of $PRP$ and $PBP$ would all be of type $\N$. These include the positions:
\[
RPP, \, RRPB, \, RPB^2, \, PPB.
\]
However, these are {\em all}\/ the children of $RPB$, so $RPB$ must be of type $\P$, contradicting our assumption. 
\end{proof}

Finally, for this section, we consider an extended form of the monotonic sequence game when the deck is $\mathbb Q$. In this extension, a chosen card can be inserted anywhere in the board, in other words you are allowed to choose the position as well as the value of the next element to insert in the sequence. A useful model of this game is that the players alternately choose points in the open unit square (or the plane, but using the square saves paper) subject to the condition that no two chosen points can lie on a vertical or horizontal line. The game ends when there are either $a$ points such that the segments connecting them all have positive slope, or $d$ such that the segments connecting them all have negative slope. We refer to such sequences of points as increasing or decreasing respectively.

This extra power reduces the analysis of the game to a simple parity argument owing to the following lemma:

\begin{lem}
Let a set of fewer than $rs$ points in the open unit square be given no two of which lie on a horizontal or vertical line. If the longest increasing sequence of points has length at most $r$ and the longest decreasing sequence of points has length at most $s$ then it is possible to add an additional point without creating a sequence of $r+1$ increasing or $s+1$ decreasing points.
\end{lem}

\begin{proof}
View the points as a permutation. To avoid trivialities, suppose that there is indeed an increasing subsequence of length $r$ and a decreasing subsequence of length $s$. Under these conditions, it is well known that the permutation has a decomposition into $s$ disjoint increasing subsequences, $I_1$ through $I_s$, each of length at most $r$ which can be obtained by a simple greedy algorithm. Since the number of elements of the permutation is less than $rs$, one of these subsequences, without loss of generality $I_1$, will contain at most $r-1$ points. Now consider a decomposition of the permutation into $r$ disjoint decreasing subsequences $D_1$ through $D_r$ each of length at most $s$ (which can also be obtained by a greedy algorithm). Since for each $i$ and $j$, $|D_i \cap I_j| \leq 1$ any of the $D_i$ of size $s$ must intersect each $I_j$. However,  some $D_i$ has empty intersection with $I_1$ (since there are $r$ $D$'s and at most $r-1$ points in $I_1$). Without loss of generality, suppose it is $D_1$ and note that necessarily $|D_1| < s$.

Now return to thinking of the elements of the permutation as points in the square. It is possible to find a point $(x,y)$ whose addition to $D_1$ forms a decreasing sequence, and whose addition to $I_1$ forms an increasing sequence. Such a point can be obtained by ``connecting the dots'' for $D_1$, and connecting the ends horizontally to the sides of the square. Do likewise for $I_1$ only connect the ends vertically. The resulting two paths have a point $P$ in common. Suppose that $P$ lies in a vertical or horizontal line determined by any of the finitely many points in the set. In that case, it is possible to perturb $P$ slightly, so that this is no longer true and so that $P$'s addition to $D_1$ forms a decreasing sequence, and its addition to $I_1$ forms an increasing one, without otherwise changing $P$'s relative horizontal or vertical position with respect to the elements of the set. The point $P$  thus satisfies the lemma since its addition still permits the partitioning of the set of points into $r$ decreasing sequences of size at most $s$, and $s$ increasing sequences of size at most $r$.
\end{proof}

In terms of the extended monotonic sequence game with parameters $(a,d,{\mathbb Q})$ the lemma above implies that for the first $(a-2)(d-2)$ moves neither player can be forced to play suicidally. However, at this point, by the Erd{\H{o}}s-Szekeres theorem the next move is necessarily suicidal. Since the parity of $ad$ is the same as that of $(a-2)(d-2)$ we obtain:

\begin{thm}
The extended monotonic sequence game with parameters $(a, d, {\mathbb Q})$ has type $\N$ if $ad$ is odd, and type $\P$ if $ad$ is even.
\end{thm}

\section{Observations and open problems}

It appears that the monotonic sequence game, particularly with normal termination criteria, has a fairly strong bias towards the first player. Specifically, our computational results suggest the following pair of conjectures:
\begin{itemize}
\item
For any $a \geq d \geq 3$ and all sufficiently large $n$, $\Wnor(a,d,n) = \Wmis(a,d,n) = \N$.
\item
For any $a \geq d \geq 3$, $\Wnor(a,d, {\mathbb Q}) = \N$.
\end{itemize}
We would be surprised (assuming the correctness of these conjectures) if similar results did not also hold for other infinite linear orders (not models of the theory of almost all finite partial orders) such as $\mathbb N$ or $\mathbb Z$.

In the finite form of the game it appears that the last $\D$ occurring in the sequence $\Wmis(a,d,n)$ is generally closer to position $a + d$ than to position $(a-1)(d-1)$. It would be of interest to determine a good upper bound for the position of this last $\D$ (the same observation and question applies to the sequence $\Wnor(a,d,n)$ though the computational evidence is less compelling). Likewise, the ``long run behaviour'' of these games seems to become established well before the bound obtained using the argument of Proposition \ref{prop:equivalence}. That the trailing sequences of $\N$'s observed in the computational results do generally represent long run behaviour is supported by a more detailed examination of these positions showing that there is a large central block of equivalent moves, which extends by a single element each time the deck size is increased (extensions to {\sc CGSuite} \cite{cgsuite} were used for some of these computations).

Another area of interest to investigate would be the behaviour of the extended form of the game played with a finite deck. In this form, players take turn naming pairs $(i, \pi_i)$ subject to the constraint that the chosen values form part of the graph of some permutation of $\{1,2,\ldots, n\}$ (and with termination based on increasing or decreasing sequences as normally). An equivalent formulation has the players placing non-attacking rooks on a (generalised) chessboard.

\bibliographystyle{plain}
\bibliography{monseq}

\end{document}